\documentclass[12pt]{article}
\makeatletter
\def\section{\@startsection{section}{1}{\z@}{-3.5ex plus -1ex minus -2.ex}
{2.3ex plus .2ex}{\Large\bf}}

\def\subsection{\@startsection{subsection}{2}{\z@}{-3.25ex plus
 -1ex minus -2.ex}
{1.5ex plus .2ex}{\bf}}

\def\vsn{\vskip 1pc \noindent}
\def\f{\newline}
\def\e{\varepsilon}

\def\comp{{\rm comp}}

\def\rr{{\bf R}}
\def\nn{{\bf N}}

\def\va{\vec{\alpha}}

\newcommand{\be} {\begin{equation}}
\newcommand{\ee} {\end{equation}}
\newcommand{\bd} {\begin{displaymath}}
\newcommand{\ed} {\end{displaymath}}

\newcommand{\bq}{\begin{eqnarray}}
\newcommand{\eq}{\end{eqnarray}}
\newcommand{\bqn}{\begin{eqnarray*}}
\newcommand{\eqn}{\end{eqnarray*}}
\newcommand{\ba}[1]{\begin{array}{#1}}
\newcommand{\eqa}{\end{array}}
\def\qed{
   \\[-4ex]
  \hbox to \hsize{\hfill \vrule height 1.6ex width 1.5ex
  depth -.1ex}}

\begin{document}

\bibliographystyle{alpha}

\begin{center} {\Large {\bf
Asymptotically tight worst case complexity bounds for initial-value problems with nonadaptive information
 }\footnotemark[1] } 
\end{center}
\footnotetext[1]{ ~\noindent This research was partly  
supported by the Polish
NCN grant - decision No. DEC-2017/25/B/ST1/00945   and  the Polish Ministry of Science 
and Higher Education. \vsn}               
  
\medskip
\begin{center}
{\large {\bf Boles\l aw Kacewicz \footnotemark[2]  }}
\end{center}
\footnotetext[2]{ 
\begin{minipage}[t]{16cm} 
 \noindent
{\it AGH University of Science 
and Technology, Faculty of Applied Mathematics,\\
\noindent  Al. Mickiewicza 30, paw. A3/A4, III p., 
pok. 301,\\
 30-059 Krakow, Poland 
\newline
E-mail B. Kacewicz: $\;\;$ kacewicz@agh.edu.pl }  
\end{minipage} }

\thispagestyle{empty}
$~$
\begin{center}  28.01.2018 \end{center}
\vsn
\vsn
\begin{center} {\bf{\Large Abstract}} \end{center}
\noindent
It is known  that, for systems of initial-value problems, algorithms using adaptive information
       perform much better in the worst case setting than  the algorithms using nonadaptive information.
In the latter case,  lower  and upper complexity bounds significantly depend on the number of equations.  
However, in contrast with adaptive information, existing lower and upper complexity bounds  
 for nonadaptive information are not asymptotically tight.  In this paper, we close the gap in the complexity 
exponents,  showing 
asymptotically matching bounds for nonadaptive standard information, as well as for a more general class
of nonadaptive linear information. 
\newpage
\noindent
\section {\bf{\Large Introduction}} 
\noindent
We aim at closing a gap between upper and lower worst case complexity bounds for initial-value problems 
with nonadaptive information. A motivation comes from a discussion on this subject that we had  with Stefan Heinrich in 2016.
We deal with the solution of systems  
\be
z'(t)=f(z(t)),\;\; t\in [a,b],\;\;\;\; z(a)=\eta,
\label{1}
\ee
where $a<b$, $f:\rr^d \to \rr^d$  is a $C^{r}$ function,  $d\geq 1$,  and $\eta\in \rr^d$.
A class of functions $f$,  denoted by $F_{r,d}$, is given by  (\ref{class}).  
For $\e>0$, we measure the difficulty of the problem  by 
the minimal 
cost of an algorithm, based on some information, that gives us an $\e$-approximation to the solution
 (the $\e$-complexity of the problem). If adaptive information is allowed, then  the $\e$-complexity is denoted by
 $\comp(\e,F_{r,d})$.  
 The notation     $\comp^{{\rm nonad}}(\e,F_{r,d})$, where  the superscript is added,  means that we restrict ourselves to 
the class of nonadaptive information. 
For precise definitions of basic notions, the reader is referred to the next section. 
Our aim is to establish
the asymptotics of $\comp^{{\rm nonad}}(\e,F_{r,d})$ as $\e\to 0$ for nonadaptive information, as function of the regularity $r$ and the dimension $d$.
\f
A question about potential advantages of adaptive over nonadaptive algorithms for solving various problems is 
an important issue in numerical analysis. Many different points of view cause some discussions  and sometimes misunderstandings
among numerical analysts in that respect.  From practical point of view, adaption is claimed to be definitely  better,  which is supported by
results of numerical experiments, see e.g. \cite{jackiewicz}, \cite{mazzia}, \cite{Pie} and many other papers. A closer look however shows
that advantages of adaption depend very much on the problem itself and the  class of problem instances being solved.
 It is not a purpose of this paper to discuss the adaption/nonadaption issue in details --
 to have  a flavor of  it, one can consult the monograph \cite{TWW}, or recent papers \cite{hickernell}, 
\cite{NumAlg}, \cite{plaskota}.
\f
In what follows, for a positive function $\gamma=\gamma(\e)$, the asymptotic expressions  $O(\gamma(\e))$,  $\Omega(\gamma(\e))$
and $\Theta(\gamma(\e))$ will always be  meant  as $\e\to 0$.
It is  known for many years that  for problem  (\ref{1}) 
adaptive information is much more efficient in the worst case setting than nonadaptive one.
It was shown for adaptive information that
the $\e$-complexity of (\ref{1}) is,  (see \cite{Kac1}):
\vsn
\begin{tabular}{ll}
$\comp(\e,F_{r,d})= \Theta\left((1/\e)^{1/r} \right),$ & \mbox{ for the class of standard adaptive information, }\\
$\comp(\e,F_{r,d})=\Theta\left((1/\e)^{1/(r+1)} \right),$ & \mbox{ for the class of linear adaptive information. }
\end{tabular}
\vsn
 In both cases of standard and linear information, the complexity bounds are  asymptotically tight, and 
the asymptotics is independent of $d$.
\f
In the nonadaptive case, the existing complexity bounds are not tight.  In \cite{Kac2}, we considered the class $F_{r,d}$ with 
$M=(0,1)^d$ and $D=1$, see (\ref{class}).
 It was shown (translating the results from non-autonomuos problems in \cite{Kac2} to the autonomous 
ones (\ref{1})) that 
\vsn
\begin{tabular}{ll}
$(a)\;\;\comp^{{\rm nonad}}(\e,F_{r,d})= \Omega\left(  (1/\e)^{d/(r+1)}   \right),$  &  \mbox{ for the class of all linear nonadaptive information, }\\
$(b)\;\; \comp^{{\rm nonad}}(\e,F_{r,d})= O\left(  (1/\e)^{d/r}   \right),$  &  \mbox{ for the class of  standard  nonadaptive information. }\\
                                                                                                                   
\end{tabular}
\vsn
The influence of the dimension  $d$ in the nonadaptive case is very significant, which  indicates that the problem (\ref{1}) 
is not well suited for  nonadaptive solution. The complexity
radically increases (asymptotically) with~$d$.
\f
The bounds {\it (a)} and {\it (b)} do not match, so that the asymptotics of the $\e$-complexity of the nonadaptive solution of
(\ref{1}) is not known. 
 In this paper,  we close the gap between lower and upper bounds in  some important cases. We  show that for $d\geq 2$
\vsn
\begin{tabular}{ll}
$\comp^{{\rm nonad}}(\e,F_{r,d})=\Omega\left(  (1/\e)^{(d-1)/r}   \right), $  &  \mbox{ for the class of all linear nonadaptive information. }
\end{tabular}
\vsn
For $d>r+1$ this is an improvement over the lower bound {\it (a)}. It shows 
in particular that it is not possible in general, as one may expect,
to achieve the complexity proportional to $(1/\e)^{d/(r+1)}$ by allowing nonadaptive linear (nonstandard) information.
Our main result, contained in Theorem 1 and next extended  in Theorem 2,  states that
\vsn
\begin{tabular}{ll}
$\comp^{{\rm nonad}}(\e,F_{r,d})=\Omega\left((1/\e)^{d/r} \right),$ & \mbox{ for a class of linear nonadaptive information }\\
     $~$                                                                               & \mbox{ that includes any standard information. }
\end{tabular}
\vsn
This  improves the lower bound {\it (a)}, and matches the upper  bound {\it (b)}.
The question about the asymptotics of the $\e$-complexity for the class of all linear nonadaptive information is still  open.
It is thought to be $\Theta\left((1/\e)^{d/r} \right)$ as $\e\to 0$, the same as for standard information.
\f
Finally, in  Remark 1 we point out how the proof of Theorem 1 can be modified to get  the   complexity lower bound for non-autonomous systems.
In the lower bound of Theorem 1, one  needs to replace $d$ by $d+1$ in the exponent. 
\f
The paper is organized as follows. In Section 2 basic notation is established and known results are recalled. Section 3 is devoted
to  the proof of  the main result  in Theorem 1, and to its extension in Theorem 2.  
In the Appendix we give auxilliary constructions and bounds.
\vsn
\section {\bf{\Large  Preliminaries}} 
\subsection {\bf{ Basic definitions}} 
\noindent
Let $r\in \nn$ ($r\geq 1$), $D>0$ and $M$ be a nonempty  open subset of $\rr^d$. We consider the class of functions
 $f$ with $r$ continuous bounded derivatives
$$
F_{r,d}=\{f: \rr^d\to \rr^d:\; f\in C^r(\rr^d), \; f(y)=0 \mbox{ for } y\notin M,\; 
$$
\be
\;\;\;\; |D^if^j(y)|\leq D, \; y\in \rr^d,\; i=0,1,\ldots,r, \; j=1,2,\ldots,d \},
\label{class}
\ee
where $ D^if^j$ denotes any partial derivative of order $i$ of $j$th component of $f$, $f=[f^1,f^2,\ldots,f^d]^T$. 
In particular, $f$ is a Lipschitz function in $\rr^d$, with the Lipschitz constant denoted by $L$, $L=L(D)$. 
We assume that $\eta\in M$, which is the only case of interest (otherwise, $z(t)\equiv \eta,\; t\in [a,b]$, for any $f\in F_{r,d}$). 
We shall use in what follows the maximum norm $\|\cdot\|$  in $\rr^d$.
\f
Let $n\in \nn$. The function $f$ is accessible through information given by a vector 
\be
N_n(f)=[L_1(f), L_2(f),\ldots, L_n(f)]^T,
\label{information}
\ee
where $L_j$ are linear functionals on $C^r(\rr^d)$. We will be interested in the power of nonadaptive information, which is defined by
functionals $L_j$  selected  simultaneously  in advance, before computation starts.  Otherwise, 
if successive functionals are selected depending on previously computed values, the information is called adaptive.
Most often, the functionals are defined by the values of $f$ or its
partial derivatives evaluated at certain points,
\be
L_j(f)= D^{k_j} f^{i_j}(y_j),\;\; \; j=1,2,\ldots,n, 
\label{standard}
\ee
where $ D^{k_j}$ is some partial derivative of order $k_j$, $0\leq k_j\leq r$, $1\leq i_j\leq d$, and $y_j\in \rr^d$. Such information is called standard. It is nonadaptive if
the points $y_j$,  as well as $k_j$ and $i_j$, $j=1,2,\ldots, n$, are given in advance.
\f
By an algorithm we mean any function $\Phi_n$ that computes 
an approximation $l:[a,b]\to \rr^d$  to the solution $z$,  based on $N_n(f)$, $l(t)= \Phi_n(N_n(f))(t)$. 
The worst case error of $\Phi_n$ with information $N_n$ in the class $F_{r,d}$ is defined by
\be
e(\Phi_n, N_n, F_{r,d})= \sup\limits_{f\in F_{r,d}}\sup\limits_{t\in [a,b]}\, \| z(t)- l(t)\|.
\label{error}
\ee
For $\e>0$, by the $\e$-(information) complexity of the problem (\ref{1}) we mean the minimal number of functionals sufficient for 
approximating $z$ with error at most $\e$, that is,
\be
\comp(\e, F_{r,d}) =\min \left\{ n\geq 1: \mbox{ there exist } N_n \mbox{ and } \Phi_n \mbox{ such that }  e(\Phi_n,N_n, F_{r,d})\leq \e\right\},
\label{comp}
\ee
with $\min \emptyset = +\infty$. 
The complexity obviously depends on the class of admitted information operators. 
If we restrict ourselves to the class of nonadaptive information, then the $\e$-complexity
will be denoted by   $\comp^{{\rm nonad}}(\e, F_{r,d})$. 
An interesting question,
in particular, is whether adaption has advantages over nonadaption.
\vsn
\subsection {\bf{  Known results on adaption versus nonadaption for IVPs }} 
\noindent
We briefly recall what is known about the potential of nonadaption for (\ref{1}).  It has been shown in 
 \cite{Kac1} and  \cite{Kac2}   that the following holds. 
 A straightforward modification is only needed
in the nonadaptive case with respect to \cite{Kac2}, where nonautonomous problems with nonadaptive infomation were considered. 
\vsn
{\bf Theorem A} (\cite{Kac1}, \cite{Kac2}) $\;\;${\it For the class of standard adaptive information
\be
\comp(\e, F_{r,d})=\Theta\left(\left(1/\e\right)^{1/r}\right).
\label{A1}
\ee
For the class of linear adaptive information
\be
\comp(\e, F_{r,d})=\Theta\left(\left(1/\e\right)^{1/(r+1)}\right).
\label{A11}
\ee
For the class of linear nonadaptive information, 
for the function class $F_{r,d}$ with $M=(0,1)^d$ and $D=1$,  we have 
\be
\Omega\left(\left(1/\e\right)^{d/(r+1)}\right) =  \comp^{{\rm nonad}}(\e, F_{r,d}) = O \left(\left(1/\e\right)^{d/r}\right).
\label{A2}
\ee
The upper bound in (\ref{A2}) is achieved by  standard information.
}
\vsn
The bounds show  a substantial advantage of adaptive over nonadaptive information for the problem (\ref{1}):
the number of equations $d$  significantly increases the complexity, if we restrict ourselves to nonadaptive information.
\f
Because of the gap between lower and upper bounds in ({\ref{A2}),  the asymptotics of the 
complexity in the nonadaptive case as $\e\to 0$  is not known (for both standard and linear information). 
We shall remove the gap in the next section for standard information,  as well as  for a class of linear information. 
\vsn
\section {\bf{\Large  Asymptotically tight complexity bounds }} 
\subsection {\bf  Standard information } 
\noindent
In the following theorem we prove new lower complexity bound improving that in  (\ref{A2}) for standard information. 
\vsn
{\bf Theorem 1}$\;\;$ {\it 
For the class of standard nonadaptive information
\be
 \comp^{{\rm nonad}}(\e, F_{r,d})= \Omega\left( \left(1/\e\right)^{d/r}\right).
\label{th1}
\ee
}
{\bf Proof}$\;$
In view of (\ref{A1}), it is enough  to consider $d\geq 2$.
It suffices to show that there exists a positive constant $C_3$ such that for sufficiently large $n$,
for any standard nonadaptive information $N_n$ and any algorithm $\Phi_n$,
 we have 
\be
e(\Phi_n,N_n, F_{r,d})\geq C_3 n^{-r/d}.
\label{12}
\ee
It is easy to see that for any $f_1,f_2\in F_{r,d}$ such that $N_n(f_1)=N_n(f_2)$, for the respective solutions $z_1$ and $z_2$
of (\ref{1}) it holds
\be
e(\Phi_n,N_n, F_{r,d})\geq (1/2) \sup\limits_{t\in [a,b]} \|z_1(t)-z_2(t)\|.
\label{13}
\ee
We shall now construct funtions $f_1$ and $f_2$ to get the lower bound in (\ref{13}) as large as possible.
Let $M_1$ be an open set containing $\eta$, whose closure is contained in $M$. We define $f_1$ such that  
\be
f_1(y)=\left\{
\begin{array}{ll} \va &\;\; \mbox{ for } y\in M_1,\\
                              0 &\;\; \mbox{ for } y\notin M.
\end{array}
\right.
\label{f1}
\ee
Here $\va$ is a vector in $\rr^d$ with  $\| \va\|=\Delta>0$ (the vector with components $[\alpha^1,\ldots, \alpha^d]^T$
may be identified with a point with the same components). The number $\Delta$ is sufficiently small to assure that:
\f
1. $f_1$ can be extended to $M \setminus M_1$  so that    $f_1\in F_{r,d}$, 
\f
2. $\Delta\leq D/2$,
\f
3. $B_\Delta = \{y:\, \|y-\eta\|\leq \Delta (b-a)\} \subset M_1$, which assures that the solution $z_1$ of (\ref{1}) for $f_1$ is given by
$$
z_1(t)=\va (t-a)+\eta, \;\; t\in [a,b].
$$
The direction of $\va$ will be chosen later on.
\f
The function $f_2$ will be given as $f_2=f_1+H$, where $N_n(H)=0$.  We now construct   $H$. 
Let $0<T\leq \Delta (b-a)$. Consider a hypercube $K$ contained in $ B_\Delta$,
\be
K=[a_1,b_1]\times [a_2,b_2]\times \ldots \times  [a_d,b_d]
\label{cubeK}
\ee
such that $[b_1,b_2,\ldots, b_d]^T=\eta$, and the length of the egdes is $T$, $b_j-a_j=T$, $j=1,2,\ldots,d$.
\f
Let $m=\lceil n^{1/d} \rceil$  and $p\in \nn$. We divide the egdes $[a_j,b_j]$, $j=1,2,\ldots,d-1$
into $pm$ intervals using equidistant points $t_k^j= a_j+ kT/(pm)$, $k=0,1,\ldots,pm$. The hypercube in 
the $(d-1)$-dimensional hyperplane 
$$\{y=[y^1,y^2,\ldots, y^d]^T:\;\; y^d=a_d \}$$
given by
 $$
K_1=[a_1,b_1]\times [a_2,b_2]\times \ldots \times  [a_{d-1},b_{d-1}]\times [a_d,a_d]
$$
is thus divided into $(pm)^{d-1}$ smaller  hypercubes 
$$
K_{k_1,k_2,\ldots, k_{d-1}} =[t_{k_1}^1, t_{k_1+1}^1]\times [t_{k_2}^2, t_{k_2+1}^2]\times\ldots \times  
[t_{k_{d-1}}^{d-1}, t_{k_{d-1}+1}^{d-1}]\times [a_d,a_d],
$$
 $k_j=0,1,\ldots, pm-1$.
\f
 Let  $c_{k_1,k_2,\ldots, k_{d-1}}\in \rr^d$ be the center of $K_{k_1,k_2,\ldots, k_{d-1}}$. 
The direction of  $\va$ in (\ref{f1}) will be selected
from among $(pm)^{d-1}$ directions of the vectors  
$$\va_{k_1,k_2,\ldots, k_{d-1}} = c_{k_1,k_2,\ldots, k_{d-1}} -\eta,\;\;\; \|\va_{k_1,k_2,\ldots, k_{d-1}}\|=T. $$
We now associate with each $\va_{k_1,k_2,\ldots, k_{d-1}}$ a parallelepiped $P_{k_1,k_2,\ldots, k_{d-1}} \subset \rr^d$ as follows.
Let $C$ be a cone with vertex $\eta$ defined as a convex hull of $K_{k_1,k_2,\ldots, k_{d-1}}$ and $\eta$.  Denote by 
 $\bar K_{k_1,k_2,\ldots, k_{d-1}}$  the intersection of $C$ with the $(d-1)$-dimensional hyperplane $y^d=(b_d+a_d)/2$.
We define 
\be
P_{k_1,k_2,\ldots, k_{d-1}} = \left\{ y\in \rr^d:\; y=\bar y  + (1/2)    \va_{k_1,k_2,\ldots, k_{d-1}}\, l ,  \mbox{ for some }
\bar y\in       \bar K_{k_1,k_2,\ldots, k_{d-1}} \mbox{ and } l\in [0,1]         \right\} .
\label{14}
\ee
Note that the paralellepipeds $P_{k_1,k_2,\ldots, k_{d-1}}$ are contained in $B_\Delta$, and have disjoint interiors for different
vectors
$(k_1,k_2,\ldots, k_{d-1})$. We finally divide each $P_{k_1,k_2,\ldots, k_{d-1}}$ into $pm$ smaller parallelepipeds (cells)
with disjoint interiors 
$$
P^0_{k_1,k_2,\ldots, k_{d-1}}, \; P^1_{k_1,k_2,\ldots, k_{d-1}}, \ldots, P^{pm-1}_{k_1,k_2,\ldots, k_{d-1}},
$$
by setting
$$
P^j_{k_1,k_2,\ldots, k_{d-1}}=\left\{ y\in P_{k_1,k_2,\ldots, k_{d-1}}:\;  l\in \left[ j/(pm), (j+1)/(pm) \right] \right\},
$$
$j=0,1,\ldots, pm-1.$ The total number of cells is $(pm)^d$.
\f
We now apply in each  $P^j_{k_1,k_2,\ldots, k_{d-1}}$  the construction described in point I of the Appendix  to get 
bump functions $\hat H^j_{k_1,k_2,\ldots, k_{d-1}}:\rr^d\to \rr$ with  supports $P^j_{k_1,k_2,\ldots, k_{d-1}}$. 
We use (\ref{bump}) with $P=P^j_{k_1,k_2,\ldots, k_{d-1}}$, $\; q_s-p_s= T/(2pm)$, $s=1,2,\ldots, d-1$,    and $\; N=T/(2pm)$.
\f
The function $H:\rr^d\to \rr^d$ is defined by
\be
H(y)=\left[\hat H(y),0,\ldots,0\right]^T,\;\;\; y\in \rr^d,
\label{funH}
\ee
where 
\be
\hat H(y)= \sum\limits_{k_1,k_2,\ldots, k_{d-1}, j} \beta_{k_1,k_2,\ldots, k_{d-1}}^j \hat H^j_{k_1,k_2,\ldots, k_{d-1}}(y) 
\label{funhatH}
\ee
with $|\beta_{k_1,k_2,\ldots, k_{d-1}}^ j|\leq 1$. The support of 
the $C^\infty$ function $H$ is the sum of all parallelepipeds (cells) $P^j_{k_1,k_2,\ldots, k_{d-1}}$.
The number of parameters $\beta_{k_1,k_2,\ldots, k_{d-1}}^j$ is $(pm)^d$.
\f
Coefficients $\beta_{k_1,k_2,\ldots, k_{d-1}}^j$ are now selected to assure that $N_n(H)=0$, where $N_n$ is  standard
information given by (\ref{information}) and (\ref{standard}). It suffices to eliminate all cells containing  in the interior
any of the information points $y_j$, that is, to set all corresponding $\beta$'s to $0$. The remaining $\beta$'s, 
whose number is at least $(pm)^d-n$, are set to $1$. 
\f
 Since $(pm)^d-n\geq (pm)^d-m^d=(p^d-1)m^d$, one can observe for $p\geq 2$, 
that there must exist $(k_1,k_2,\ldots, k_{d-1})$ such that 
\be
\beta_{k_1,k_2,\ldots, k_{d-1}}^j=1 \mbox{  for at least } 1/2 \mbox{ of indices } j=0,1,\ldots,pm-1. 
\label{beta}
\ee
We now select $\va$ in (\ref{f1}) to be 
\be
\va=  \frac{\va_{k_1,k_2,\ldots, k_{d-1}} }{\| \va_{k_1,k_2,\ldots, k_{d-1}}  \| }\, \Delta,   
\label{122}
\ee
for the chosen indices $k_1,k_2,\ldots, k_{d-1}$. 
\f
 Summarizing the above construction, 
functions $f_1$ and $f_2=f_1+H$ belong to  $F_{r,d}$,  $N_n(f_1)=N_n(f_2)$, and the solution 
for $f_1$ has the form $z_1(t)=\va (t-a)+\eta, \; t\in [a,b].$ 
To complete the proof, it remains to bound from below the distance $\sup\limits_{t\in [a,b]} \|z_1(t)-z_2(t)\| $. From (\ref{lowbound}),
$$
\sup\limits_{t\in [a,b]} \|z_1(t)-z_2(t)\| \geq \frac{1}{1+L(b-a)} \left\| \int\limits_a^b H(\va(\xi-a)+\eta)\, d\xi\right\| 
$$
\be
= \frac{1}{1+L(b-a)} \left| \int\limits_a^b \hat H(\va(\xi-a)+\eta)\, d\xi\right| .
\label{111}
\ee
The values  $\hat H(\va(\xi-a)+\eta)$  are $0$ if the argument is outside the parallelepiped   $P_{k_1,k_2,\ldots, k_{d-1}}$. 
The line  $\va(\xi-a)+\eta$ intersects the hyperplanes $y^d=(a_d+b_d)/2$ and $y^d=a_d$  for
$\xi=t_1=a+T/(2\Delta)$ and $\xi=t_2=a+T/\Delta \, (\leq b)$, respectively.  Hence,
$$
 \int\limits_a^b \hat H(\va(\xi-a)+\eta)\, d\xi = \int\limits_{t_1}^{t_2} \hat H(\va(\xi-a)+\eta)\, d\xi
$$
\be
=\sum\limits_{j=0}^{pm-1}     \int\limits_{\xi_j}^{\xi_{j+1}} \hat H(\va(\xi-a)+\eta)\, d\xi, 
\label{sumacalek}
\ee
where $\xi_j=t_1+ (t_2-t_1)j/(pm)= a+T/(2\Delta) + jT/(2\Delta pm)$, $j=0,1,\ldots, pm$. 
\f
Note that for $\xi\in [\xi_j, \xi_{j+1}]$, and only for such,  the argument $\va(\xi-a)+\eta$ remains in $P_{k_1,\ldots, k_{d-1}}^j$.
\f
We now use (\ref{bump}). By construction, in any cell for which $\beta_{k_1,k_2,\ldots, k_{d-1}}^j=1$, we have for $m\geq T/(2p)$
\be
\hat H(\va(\xi-a)+\eta) = \hat C \left( T/(2pm)\right)^r \left(h(1/2)\right)^{d-1} h(\bar l/N),\;\;\; \xi\in [\xi_j,\xi_{j+1}], 
\label{hatH}
\ee
where $N=T/(2pm)$ and  $\bar l/N= (\xi-\xi_j)/(\xi_{j+1}-\xi_j)$. 
Hence, we have that
$$
\sum\limits_{j=0}^{pm-1}     \int\limits_{\xi_j}^{\xi_{j+1}} \hat H(\va(\xi-a)+\eta)\, d\xi = 
 \hat C \left( T/(2pm)\right)^r \left(h(1/2)\right)^{d-1} \sum_{j-}  \int\limits_{\xi_j}^{\xi_{j+1}} h( (\xi-\xi_j)/(\xi_{j+1}-\xi_j) )\,d\xi
$$
\be
= \hat C \left( T/(2pm)\right)^r \left(h(1/2)\right)^{d-1} \int\limits_0^1 h(x)\, dx \, \sum_{j-}  (\xi_{j+1}-\xi_j),
\label{sumacalek1}
\ee
where the sum $\sum_{j-}$ extends over all indices $j$ for which $\beta_{k_1,k_2,\ldots, k_{d-1}}^j=1$. By construction, 
there is at least $pm/2$ such indices $j$, so that  the sum $\sum_{j-}  (\xi_{j+1}-\xi_j)$ is a positive constant.    This allows us to conclude that
$$
\left| \int\limits_a^b \hat H(\va(\xi-a)+\eta)\, d\xi\right| =\Omega\left(m^{-r}\right).
$$
In view of (\ref{111}) and the fact that $m=\Theta(n^{1/d})$, we get (\ref{12}),  which completes the proof. \qed
\vsn
\subsection {\bf{  Generalization to a class of linear information }} 
\noindent
Let a sequence $\{k(n)\}_{n=1}^\infty \subset \nn$   be such  that there exist $\alpha>0$ and $n_0\in \nn$ such that
\be
k(n)\leq \alpha n^{1-1/d}\;\;\; \mbox{ for } n\geq n_0.
\label{181}
\ee
We consider the following class of linear information operators
\be
N_n(f)=[\bar L_1(f),\ldots, \bar L_{n-k}(f), L_1(f),\ldots, L_k(f)]^T, \;\;\; n\geq k+1, 
\label{liniowagen}
\ee
where $k=k(n)$,  and $\bar L_j$ are  arbitrary standard information functionals defined by (\ref{standard}).  
The functionals  $L_j$ are arbitrary linear. 
\f  
In the following theorem we first somewhat improve the lower bound in (\ref{A2}) 
of Theorem~A in the class of all linear nonadaptive information. The bound  (\ref{th2a}) will yield that   we cannot in general
reduce the exponent in the upper complexity bound from $d/r$ to $d/(r+1)$,  when switching from nonadaptive standard to 
nonadaptive linear  information (as one may expect).
Next, more interesting,  
we generalize  Theorem~1 to the class of information given by  (\ref{liniowagen}),   
by proving  that  the lower bound  (\ref{th1}) still
holds true in this case. 
\vsn
{\bf Theorem 2}$\;\;$ {\it  (i)$\;\;$ For the class of linear nonadaptive information (\ref{information}), we have 
\be
 \comp^{{\rm nonad}}(\e, F_{r,d})= \Omega\left( \left(1/\e\right)^{\max\{(d-1)/r, d/(r+1)\} }\right).
\label{th2a}
\ee
\noindent
(ii)$\;\;$  Let  $\{k(n)\}_{n=1}^\infty$ satisfy (\ref{181}).  For the class of  nonadaptive information (\ref{liniowagen}),    we have
\be
 \comp^{{\rm nonad}}(\e, F_{r,d})= \Omega\left( \left(1/\e\right)^{d/r}\right).
\label{th2}
\ee
(The constant and the maximal value of $\e$ in the ' $\Omega$' notation in (\ref{th2})  depend on $\alpha$ and $n_0$.)}
\f
{\bf Proof } of {\it (i)} $\;\;$  
In view of Theorem A, it suffices to show that for $d\geq 2$ we have 
$$
\comp^{{\rm nonad}}(\e, F_{r,d})= \Omega\left( \left(1/\e\right)^{(d-1)/r }\right).
$$
We refer to the proof of Theorem 1.  We now choose $m=\lceil (n+1)^{1/(d-1)} \rceil$ and $p=1$, and
construct the parallelepipeds $ P_{k_1,\ldots, k_{d-1}}$ as in (\ref{14}). We do not divide them further into small cells, and
define the function $H$ by (\ref{funH}) with $\hat H$ given by
\be
\hat H(y)= \sum\limits_{k_1,k_2,\ldots, k_{d-1}} \beta_{k_1,k_2,\ldots, k_{d-1}} \hat H_{k_1,k_2,\ldots, k_{d-1}}(y) 
\label{funhatH1}
\ee
with $|\beta_{k_1,k_2,\ldots, k_{d-1}}|\leq 1$.  
Here the functions $\hat H_{k_1,k_2,\ldots, k_{d-1}}$ are defined in a similar way as $\hat H^j_{k_1,k_2,\ldots, k_{d-1}}$
in the proof of Theorem 1 (now with supports $ P_{k_1,\ldots, k_{d-1}}$).
\f
The number of unknowns  $\beta_{k_1,k_2,\ldots, k_{d-1}}$ is now
$m^{d-1}\geq n+1$. The condition $N_n(H)=0$ leads to $n$ linear homogenous equations, which have 
a solution with the maximum modulus coefficient $\beta_{k_1,k_2,\ldots, k_{d-1}}=1$ for some $k_1,k_2,\ldots, k_{d-1}$.
We use these indices in (\ref{122}), and follow the further steps of the proof of Theorem 1. The sum of the integrals in
(\ref{sumacalek})  contains now only one element.
We finally get that
$e(\Phi_n, N_n, F_{r,d})=\Omega\left( m^{-r}\right)$ for any information $N_n$ and algorithm $\Phi_n$, which means that
\be
e(\Phi_n,N_n, F_{r,d})=\Omega\left( n^{-r/(d-1)}\right), \; n\to\infty. 
\label{129}
\ee
This leads to the desired bound on the complexity.  
\vsn
{\bf Proof } of {\it (ii)} $\;\;$ Let $d\geq 2$.    We repeat the steps of the proof  of Theorem 1 up to 
the definition of the function  $H$ in (\ref{funH}).  
A difference is in the definition of the coefficients $\beta_{k_1,k_2,\ldots, k_{d-1}}^j$ in order to assure that $N_n(H)=0$.
We remember that 
each of  $\beta_{k_1,k_2,\ldots, k_{d-1}}^j$ is related to a cell $P^j_{k_1,k_2,\ldots, k_{d-1}}$, which is a part 
of a parallelepiped  $P_{k_1,k_2,\ldots, k_{d-1}}$ (number these parallelepipeds is $(pm)^{d-1}$).
 To fulfill the standard information conditions $\bar L_j(H)=0$, $j=1,2,\ldots, n-k$, it suffices to eliminate cells which contain 
in the interior any of the information points that define the functionals $\bar L_j$. This is done by setting
the corresponding $\beta$'s to $0$. The number of eliminated cells is at most $n-k$. We also exclude all 
parallelepipeds $P_{k_1,k_2,\ldots, k_{d-1}}$  which  contain more than $pm/2$ eliminated cells. Denoting the number 
of these by $x$, we observe that $x\leq 2(n-k)/(pm)$. Indeed,  the number of eliminated cells is at least $pmx/2$.
If $x> 2(n-k)/(pm)$, then the number of eliminated cells is greater than $n-k$, which is a contradiction.
\f
The number of remaining parallelepipeds is thus at least   $s=(pm)^{d-1} - 2(n-k)/(pm)$. In any such parallelepiped 
$P_{k_1,k_2,\ldots, k_{d-1}}$, the number of coefficients $\beta_{k_1,k_2,\ldots, k_{d-1}}^j=0$ is at most $pm/2$.
We set the remaining $\beta_{k_1,k_2,\ldots, k_{d-1}}^j$ to be the same for each  $j$, 
$\beta_{k_1,k_2,\ldots, k_{d-1}}^j \equiv \beta_{k_1,k_2,\ldots, k_{d-1}}$.  The number of unknown coefficients  $\beta_{k_1,k_2,\ldots, k_{d-1}}$
is thus at least $s$.
\f
We select them to satisfy the conditions $L_j(H)=0$, $j=1,2,\ldots, k$.  This is the system  of $k$
linear homogenous equations with at least $(pm)^{d-1} - 2(n-k)/(pm)$ unknowns. Since $\{k(n)\}_{n=1}^\infty$ satisfies 
(\ref{181}),   we have  for sufficiently large $n$ ($m$) that 
\be
k< (pm)^{d-1} - 2(n-k)/(pm).
\label{225}
\ee
Indeed,  the condition  $k\leq \alpha n^{1-1/d}$ for sufficiently large $n$ yields that  $k\leq \alpha m^{d-1}$ for  sufficiently  large $m$.
A sufficient condition for (\ref{225}) is 
$$
\alpha m^{d-1} < \frac{(1/2) p^d -1}{ (1/2)pm-1} m^d,
$$
which is fulfilled if 
$$
\alpha \leq p^{d-1}- 2/p.
$$
The last condition holds true for sufficiently large (fixed) $p$.
\f
By (\ref{225}),  the system of equations has a  solution with the maximum modulus component $\beta_{k_1,k_2,\ldots, k_{d-1}}=1$ for some 
$k_1,k_2,\ldots, k_{d-1}$. With this choice of $k_1,k_2,\ldots, k_{d-1}$, we repeat, starting from (\ref{122}), the remaining steps 
of the proof of Theorem 1. We use the fact that $\beta^j_{k_1,k_2,\ldots, k_{d-1}}=1$ for at least $pm/2$ indices $j$,
while the remaining  $\beta^j_{k_1,k_2,\ldots, k_{d-1}}$ are $0$.
\vsn
Let $d=1$. 
We show that for each $\alpha$ there are $C>0$ and $\bar n_0$ such that for any $n\geq \bar n_0$, 
  any information $N_n$ given by (\ref{liniowagen})
 with  $k\leq \alpha$, and any algorithm $\Phi_n$ it holds
\be
e(\Phi_n,N_n, F_{r,1})\geq Cn^{-r}.
\label{71}
\ee
We only sketch the proof, by showing how to define 'difficult' functions
$f_1$ and $f_2$ in $F_{r,1}$ with the same information, $N_n(f_1)=N_n(f_2)$.  Take an interval $[\eta, \eta+\delta]$, $\delta>0$,
contained in $M$. Take $f_1$ such that  $f_1(y)\equiv \hat\alpha >0$, $y\in [\eta,\eta+\delta]$. We choose $\hat\alpha$ small enough to
assure that: $f_1$ can be extended to $\rr$ so  that the extension is in $F_{r,1}$, 
$\hat\alpha\leq D/2$ and $\hat\alpha (b-a)\leq \delta$. Then 
the solution for $f_1$ in $[a,b]$ is given by $z_1(t)=\eta + \hat\alpha(t-a),\; t\in [a,b]$. 
\f
We now divide $[\eta,\eta+ \hat\alpha (b-a)]$ uniformly into $k+1$  subintervals $I_j$, $j=1,2,\ldots, k+1$. Next, each $I_j$
is further divided into $2n$ equal subintervals $I_j^l$, $l=1,2,\ldots, 2n$.  We define function $f_2=f_1+H$
(with the solution $z_2$), where $H$ 
is given as follows. 
\f
In each subinterval $I_j^l$ containing in the interior an information point (see the functionals $\bar L_p$), we set $H\equiv 0$. There is at most $n-k$
such 'removed' subintervals $I_j^l$.  Hence, there is at most $n-k$ removed (and at least $n+k$ remaining)  subintervals in each interval $I_j$.
The conditions   $\bar L_p(H)=0$, $p=1,2,\ldots, n-k$ are automatically satisfied. 
\f
To assure that $L_j(H)=0, j=1,2,\ldots, k$, we construct in each of the remaining subintervals $I_j^l$
a standard scalar (normalized)   bump function $h_j^l$ with support $I_j^l$ (see the Appendix I), and define
$H(y)=\sum\limits_{j=1}^{k+1} \beta_j\sum\limits_l h_j^l(y)$ for $y\in \rr$.  The second sum is taken over all $l$ such that $I_j^l$ has not 
been removed. 
The conditions $L_j(H)=0$,  $j=1,2,\ldots, k$ are equivalent to $k$ linear homogeneous equations with $k+1$ unknowns $\beta_j$.
There exists a solution of the system with maximum modulus component $\beta_{ j^*}=1$.  In 
the lower bound on $\sup\limits_{\xi\in [a,b]} |z_2(\xi)-z_1(\xi)|$ given in  the Appendix II,  we choose the interval $[x,t] = I_{ j^*}$.
Since the function $H$   in the interval $I_{ j^*}$ is composed of at least $n+k$ nonzero bump functions,  the desired lower bound (\ref{71})
follows. 
\f
By inspecting the proof, we see that (\ref{71}) holds true as well for adaptive information $N_n$. 
$~$ \qed
\vsn
We finally discuss the main result of this paper  for non-autonomous systems. 
\vsn
{\bf Remark 1}  {\it  (Theorem 1 for non-autonomous problems)}   $\;$ Consider a non-autonomous system
\be
z'(t)=f(t,z(t)),\;\; t\in [a,b],\;\;\;\; z(a)=\eta,
\label{non-aut}
\ee
where  $f:[a,b]\times \rr^d \to \rr^d$  is a function from the class
$$
F_{r,d}^{{\rm non}}=\{f: [a,b]\times \rr^d\to \rr^d:\; f\in C^r([a,b]\times \rr^d), \; f(t,y)=0 \mbox{ for } y\notin M,\; 
$$
\be
\;\;\;\; |D^if^j(t,y)|\leq D, \; (t,y)\in [a,b]\times \rr^d,\; i=0,1,\ldots,r, \; j=1,2,\ldots,d \}.
\label{classnon}
\ee
We briefly sketch how  the non-autonomous case can be covered by a similar analysis as in Theorem 1,  by only 
pointing out differences  in the proof.   
The corresponding function to $f_1$ in (\ref{f1})  is now constructed  so that   
\be
f_1(t, y)=\left\{
\begin{array}{ll} \va &\;\; \mbox{ for } t\in [a,b],\;  y\in M_1,\\
                              0 &\;\; \mbox{ for } t\in [a,b],\;  y\notin M,
\end{array}
\right.
\label{f1non}
\ee
where $\|\va\| \leq \Delta$.
The construction   of the function  $H:\rr^{d+1}\to \rr^d$  corresponding to that in (\ref{funH}) 
starts in the hypercube $K\subset \rr^{d+1}$ contained in $[a,b]\times B_\Delta$ given by
\be
K=[a,b]\times [a_1,b_1]\times [a_2,b_2]\times \ldots \times  [a_d,b_d],
\label{cubeKnon}
\ee
where $B_\Delta$ is defined  in point 3.  in the proof of Theorem 1,  and $a_s$, $b_s$ are given as in (\ref{cubeK}).
Compared to the autonomous case, we have here an additional variable $t\in [a,b]$.
The hypercube $K$ will contain all graphs of the solutions  $[t, (t-a)\va^T+\eta^T]^T$, $t\in [a,b]$,   for functions $f_1$
(i.e., for vectors $\va$) under consideration. 
Let $m=\lceil n^{1/(d+1)} \rceil$.   We define the $d$-dimensional hypercube $K_1$  by
$$
K_1= [b,b]\times [a_1,b_1]\times [a_2,b_2]\times \ldots \times  [a_d,b_d].
$$
We divide $K_1$  (uniformly) into $(pm)^d$ small hypercubes $K_{k_1,k_2,\ldots,k_d}$ with centers $c_{k_1,k_2,\ldots,k_d}=
[b, \bar c_{k_1,k_2,\ldots,k_d}^T]^T $, $\bar c_{k_1,k_2,\ldots,k_d}\in \rr^d$.
The vector $\va$ will be selected as      
$$
\va= \va_{k_1,k_2,\ldots,k_d} = (\bar c_{k_1,k_2,\ldots,k_d} -\eta)/(b-a), 
$$
for a proper choice of $k_1,k_2,\ldots,k_d$.  Note that  $\|\va\| \leq \Delta$, and 
the graph $[t, (t-a)\va^T+\eta^T]^T$ passes through $[a,\eta^T]^T$ for $t=a$,
and $c_{k_1,k_2,\ldots,k_d}$ for $t=b$. 
\f
For each $d$-tuple of the indices $k_1,k_2,\ldots,k_d$, we  consider 
a cone $C\subset \rr^{d+1}$ defined as  the convex hull of $K_{k_1,k_2,\ldots,k_d}$ and $[a,\eta^T]^T$. 
 The intersection of $C$ with the $d$-dimensional hyperplane $t=(a+b)/2$ is denoted by $\bar K_{k_1,k_2,\ldots,k_d}$. 
The remaining steps  of the definition of the functions $\hat H:\rr^{d+1}\to \rr$,  $H:\rr^{d+1}\to \rr^d$,  and  
of the choice of  indices   $k_1,k_2,\ldots,k_d$ that define $\va$ are  similar to the steps 
described  in (\ref{14})--(\ref{beta}) in the 
proof of Theorem 1, with $d:=d+1$.    After replacing   $\hat H( \va(\xi-a)+\eta)$   in (\ref{111})  by
$\hat H(\xi, \va(\xi-a)+\eta) $,   the rest of the  proof goes similarly as in the autonomous case 
with $d$ replaced by $d+1$. 
\f
 Consequently,    in place of the bound in  Theorem 1, we get for non-autonomous systems  
the (matching) bound 
\be
 \comp^{{\rm nonad}}(\e, F_{r,d})= \Omega\left( \left(1/\e\right)^{(d+1)/r}\right).
\label{th1non}
\ee
\newpage\noindent
\section {\bf{\Large  Appendix}} 
\noindent
{\bf I.}$\;\;$ {\it Construction of bump functions on a parallelepiped}
\vsn
We show the construction of a $C^{\infty}$ bump function on a paralellepiped   in $\rr^d$, $d\geq 2$, which was used 
to define the function $H$ in (\ref{funH}). The construction is,
up to some details, standard. Let $h:\rr\to\rr$ be given by
\be
h(x)=\left\{ 
\begin{array}{ll} \exp\left(1/x(x-1)\right), & \;\; x\in (0,1)\\
                                                       0 & \mbox{ otherwise } .
\end{array} \right.
\label{hx}
\ee
We have that $h\in C^\infty(\rr)$, $h(x)>0$, $x\in (0,1)$, $\max\limits_{x\in \rr} h(x)=h(1/2)=c_0$,
 $\max\limits_{x\in \rr}| h^{(j)}(x)|=c_j$  for $j\geq 1$, and 
\be
\int\limits_{\rr} h(x)\, dx=\bar c,
\label{int}
\ee
where $c_j, \bar c$ are absolute positive constants. 
\f
Consider a paralellepiped  $P$ in $\rr^d$ defined as follows.
Let $P_1=[p_1,q_1]\times [p_2,q_2]\times\ldots\times [p_{d-1},q_{d-1}]$, where $q_j-p_j>0$. Let $\vec{\alpha}\in \rr^d$ be 
such that $\cos(\va, e_d)=\vec{\alpha}^T e_d/\|\va\|_2 \geq \hat c>0$, where $e_d=[0,0,\ldots,1]^T$, and $\|\cdot\|_2$ is the Euclidean norm in $\rr^d$.
For $N>0$, we set
\be
P=\left\{ y\in \rr^d:\; \mbox{ there exist } \bar y\in P_1 \mbox{ and } \bar l\in [0,N] \mbox{ such that } y=[\bar y^T,0]^T +\frac{\va}{\|\va \|} \bar l \right\}.
\label{parallel1}
\ee
We define a real-valued bump function $\hat H$ on $P$. A point $y\in \rr^d=[y^1,y^2,\ldots,y^d]^T$ is uniquely defined by the pair 
$(\bar y,\bar l)$ with $\bar y=[\bar y^1,\bar y^2,\ldots, \bar y^{d-1}]\in \rr^{d-1}$
and $\bar l\in \rr$ ($\bar y$ is a projection of $y$ on $\rr^{d-1}$ along the direction  $\va$). We set for $\hat C>0$ 
\be
\hat H(y)=\hat C \left(\min\left\{ \min\limits_{1\leq j\leq d-1}(q_j-p_j), N, 1 \right\} \right)^r \;
\prod\limits_{j=1}^{d-1} h\left((\bar y^j-p_j)/(q_j-p_j)\right) h\left(\bar l/N\right).
\label{bump}
\ee
Since the change of the variables is given by
$$
\bar y=y-\frac{\va}{\| \va \|}\cdot \frac{y^d}{ \cos(\va, e_d)}  \;\; \mbox{(the last component equal to $0$ is omitted)},\;\; 
\bar l= \frac{ y^d}{\cos(\va, e_d)},
$$ 
this describes a $C^\infty(\rr^d)$ nonnegative mapping from $\rr^d$ to $\rr$ with support $P$. Furthermore, by selecting
proper (sufficiently small, fixed) value of $\hat C$ independent of $p_j,q_j, N$, 
 we assure that $\hat H(y)\leq D/2$ and $|D^k \hat H(y)|\leq D$ 
for any partial derivative
of $\hat H$ of order $k$, $k=1,2,\ldots,r$, $y\in \rr^d$.   
\f
In the proof of Theorem 1, we use (\ref{bump}) with  $\va=-\va_{k_1,k_2,\ldots,k_{d-1}}$, and we have $\cos(\va, e_d)\geq 1/\sqrt{d}$. 
\vsn
{\bf II.}$\;\;$ {\it  A lower bound on the distance between two solutions}
\vsn
We recall for completness a standard lower bound that is used in this paper. 
Let $z_j$ be the solution of (\ref{1}) for a right-hand side $f_j$, $j=1,2 $,   $f_2=f_1+H$.
We have that for any $a\leq x\leq t\leq b$
$$
z_2(t)-z_1(t) -(z_2(x)-z_1(x))= \int\limits_x^t \left( f_2(z_2(\xi)) - f_2(z_1(\xi))\right)\, d\xi + \int\limits_x^t \left( f_2(z_1(\xi)) - f_1(z_1(\xi))\right)\, d\xi.
$$
This gives
$$
2\sup\limits_{\xi\in [x,t]} \|z_2(\xi)-z_1(\xi)\| \geq 
\left\| \int\limits_x^t H(z_1(\xi))\, d\xi \right\| - \left\| \int\limits_x^t\left(f_2(z_2(\xi)-f_2(z_1(\xi)\right)\, d\xi \right\|.
$$
From  the Lipschitz condition, we have that
\be
\sup\limits_{\xi\in [x,t]} \|z_2(\xi)-z_1(\xi)\| \geq \frac{1}{2+L(t-x)} \left\| \int\limits_x^t H(z_1(\xi))\, d\xi\right\|.
\label{lowbound}
\ee
In particular,  for any $a\leq x\leq t\leq b$, 
\be
\sup\limits_{\xi\in [a,b]} \|z_2(\xi)-z_1(\xi)\| \geq \frac{1}{2+L(t-x)} \left\| \int\limits_x^t H(z_1(\xi))\, d\xi\right\|.
\label{lowbound1}
\ee
It is easy to see that for $x=a$, the value $2+L(t-a)$ in the denominator can be replaced by $1+L(t-a)$.
\vsn
{\bf Acknowledgments}   During the visit of Stefan Heinrich to Cracow in 2016, we had a discussion on nonadaption for IVPs.
Although preliminary ideas that emerged at that time were different from those of the present  paper, the discussion 
was very much inspiring.
\vsn

\end{document}